\documentclass[12pt, reqno]{amsart}
\usepackage{graphics}
\usepackage{amsfonts}

\newcommand{\bproof}{{\bf\underline{Proof:}} }
\def\Box{\vcenter{\vbox{\hrule\hbox{\vrule
     \vbox to 8.8pt{\hbox to 10pt{}\vfill}\vrule}\hrule}}}

\newcommand{\Ff}{{\mathbb F}}
\newcommand{\Zz}{{\mathbb Z}}

\newcommand{\eproof}{\hfill$\Box$\vspace{4mm}}

\newcommand{\Z}{{\mathbb Z}}

\newtheorem{thm}{Theorem}[section]
\newtheorem{lemma}[thm]{Lemma}
\newtheorem{cor}[thm]{Corollary}

\def\Tr{\operatorname{Tr}}
\def\tr{\operatorname{tr}}
\DeclareMathOperator{\PG}{PG} \DeclareMathOperator{\GR}{GR}

\numberwithin{equation}{section}

\begin{document}

\title[Negative Latin Square Type PDS]
{Negative Latin Square Type Partial Difference Sets in
Nonelementary Abelian 2-Groups}

\author{James A. Davis, Qing Xiang}

\address{Department of Mathematics and Computer Science, University of Richmond, Richmond, VA 23173, USA,
email: {\tt jdavis@richmond.edu}}

\address{Department of Mathematical Sciences, University of Delaware, Newark, DE 19716, USA,
email: {\tt xiang@math.udel.edu}}

\insert\footins{{2000 {\it Mathematics Subject Classification} 05B10, 05E30, 11T15}}


\date{}

\begin{abstract}
Combining results on quadrics in projective geometries with an
algebraic interplay between finite fields and Galois rings, we
construct the first known family of partial difference sets with
negative Latin square type parameters in nonelementary abelian
groups, the groups $\Z_4^{2k}\times \Z_2^{4 \ell-4k}$ for all $k$
when $\ell$ is odd and for all $k < \ell$ when $\ell$ is even.
Similarly, we construct partial difference sets with Latin square
type parameters in the same groups for all $k$ when $\ell$ is even
and for all $k<\ell$ when $\ell$ is odd. These constructions
provide the first example that the non-homomorphic bijection
approach outlined by Hagita and Schmidt \cite{hagitaschmidt} can
produce difference sets in groups that previously had no known
constructions.  Computer computations indicate that the strongly
regular graphs associated to the PDSs are not isomorphic to the
known graphs, and we conjecture that the family of strongly
regular graphs will be new.
\end{abstract}
\maketitle

\section{Introduction}
\label{intro}

A $k$-element subset $D$ of a finite multiplicative group $G$ of
order $v$ is called a {\em $(v, k, \lambda, \mu)$-partial
difference set} (PDS) in $G$ provided that the multiset of
``differences'' $\{ d_{1}d_{2}^{-1} \mid d_{1}, d_{2} \in D, d_{1}
\neq d_{2} \}$ contains each nonidentity element of $D$ exactly
$\lambda$ times and each nonidentity element in $G \backslash D$
exactly $\mu$ times.  Partial difference sets are equivalent to
strongly regular graphs with a regular automorphism group, and
they are connected to projective two-weight codes and
two-intersection sets in projective spaces over finite fields. See
\cite{ma} or \cite{calderbankkantor} for background on these
alternative approaches.

Ma's survey~\cite{ma} identifies several families of PDSs. Among
these are PDSs with parameters $(n^2, r(n - \epsilon), \epsilon n
+ r^2 - 3 \epsilon r, r^2 - \epsilon r)$ for $\epsilon = \pm 1$.
When $\epsilon = 1$, the PDS is called a {\it{Latin square type
PDS}}, and when $\epsilon = -1$, the PDS is called a {\it{negative
Latin square type PDS}}.  Early results in this area can be found
in the paper by Bailey and Jungnickel~\cite{baileyjungnickel}.  We
will focus on the case where $n = 4^{2 \ell}$ and $r = 4^{\ell -
1} + \epsilon$. PDSs with these parameters have been constructed
using quadratic forms over $\Ff_4$, the field with four elements,
as we will mention in Section~\ref{quadraticforms}. This
construction, as is true of most of the constructions the authors
are aware of, are associated to elementary abelian groups.

Several authors, including those of~\cite{crx}, \cite{dx},
\cite{hou}, \cite{hlx}, \cite{leungma2}, and \cite{rx}, have used
Galois rings, and more generally finite local rings, to construct
PDSs in groups that are not elementary abelian. As far as we know,
all PDSs constructed by finite local rings except those in
\cite{dx} have Latin square parameters. After a presentation on
one of these constructions, Mikhail Klin and William Martin
\cite{klinmartin} asked whether local ring construction can
produce PDSs with negative Latin square type parameters. In this
paper, we construct the first known family of PDSs with negative
Latin square type parameters in nonelementary abelian groups by
using Galois rings, and therefore answer the aforementioned
question in the affirmative.

Recent work by Mathon~\cite{mathon} on maximal arcs in projective
planes motivated the authors to revisit their work~\cite{dx} on
Denniston parameter PDSs.  We recognized that one of the Galois
ring constructions in that paper could be obtained using a natural
map from the finite field construction. The mapping is not an
isomorphism, but it is a bijection, and there is a relationship
between the character sums in the two group rings.  We extend that
observation to higher dimensional projective spaces. Hagita and
Schmidt \cite{hagitaschmidt} recently proposed using bijections
between group rings with the property that character sums are
preserved, and that is the approach we take in this paper.
Wolfmann \cite{wolfmann} used a similar approach to constructing
Hadamard difference sets in nonelementary abelian groups via
bijections between groups: like Hagita and Schmidt's paper, the
examples found in Wolfmann's paper involved groups which were
already known to contain difference sets.  This paper uses similar
techniques to construct PDSs in groups having no known previous
constructions, the first time this approach has led to something
new.  (We note that Bruck \cite{bruck} used a bijection between
group rings to construct difference sets in nonabelian groups. His
mapping did not preserve character sums in the same way as the
Hagita-Schmidt approach.)

We will limit our attention to abelian groups in this paper.  In
that context, a (complex) character of an abelian group is a
homomorphism from the group to the multiplicative group of complex
roots of unity.  The {\em principal character} is the character
mapping every element of the group to 1.  All other characters are
called {\it nonprincipal}. Starting with the important work of
Turyn \cite{turyn}, character sums have been a powerful tool in
the study of difference sets of all types.  The following lemma
states how character sums can be used to verify that a subset of a
group is a PDS.

\begin{lemma}
\label{charsum-PDS} Let $G$ be an abelian group of order $v$ and
$D$ be a $k$-subset of $G$ such that $\{ d^{-1} \mid d \in D \} =
D$ and $1\not\in D$. Then $D$ is a $(v,k,\lambda, \mu)$-PDS in $G$
if and only if, for any complex character $\chi$ of $G$,

$$\sum_{d \in D} \chi(d) = \left\{ \begin{array}{ll}
    k & \mbox{if $\chi$ is principal on $G$.} \\
    \frac{(\lambda - \mu) \pm \sqrt{(\lambda - \mu)^2 + 4(k-\mu)}}{2} & \mbox{if $\chi$ is nonprincipal on $G$.}
            \end{array}
        \right. $$
\end{lemma}

\subsection{Projective two-intersection sets}

Let $\PG(m-1,q)$ denote the desarguesian $(m-1)$-dimensional
projective space over the finite field $\Ff_q$, where $q$ is a
power of a prime $p$, and let $\Ff_q^{m}$ be the $m$-dimensional
vector space associated with $\PG(m-1,q)$. A {\it projective}
$(n,m,h_1,h_2)$ {\it set} ${\mathcal O}=\{\langle y_1\rangle,
\langle y_2\rangle, \cdots ,\langle y_n\rangle\}$ is a proper,
non-empty set of $n$ points of the projective space $\PG(m-1,q)$
with the property that every hyperplane meets ${\mathcal O}$ in
$h_1$ or $h_2$ points. Define $\Omega =\{v\in
\Ff_q^{m}\setminus\{0\}\mid \langle v\rangle\in {\mathcal O}\}$ to
be the set of vectors in $\Ff_q^{m}$ corresponding to ${\mathcal
O}$, i.e., $\Omega=\Ff_q^*{\mathcal O}$.

For $(w_1,w_2,\ldots,w_{m}) \in \Ff_q^{m}$, we define a character
of the additive group of $\Ff_q^{m}$ as follows:
$$\chi_{_{(w_1,w_2,\ldots, w_{m})}}:
(v_1,v_2,\ldots, v_{m}) \mapsto \xi_p^{\sum_{i=1}^{m} {\rm tr}(w_i
v_i)}, (v_1,v_2,\ldots, v_{m}) \in \Ff_q^{m},$$ where $\xi_p$ is a
complex primitive $p$-th root of unity and ${\rm tr}$ is the trace
from $\Ff_q$ to $\Ff_p$. It is easy to see that
$\chi_{_{(w_1,w_2,\ldots, w_{m})}}$, $(w_1,w_2,\ldots,w_{m}) \in
\Ff_q^{m}$, are all the characters of the additive group of
$\Ff_q^{m}$.

For any nontrivial additive character $\chi_{_{(w_1,w_2,\ldots,
w_{m})}}$ of $\Ff_q^{m}$, we have
\begin{eqnarray*}\label{charcomp}
\chi_{_{(w_1,w_2,\ldots, w_{m})}}(\Omega)
&=&(q-1)|(w_1,w_2,\ldots,w_{m})^{\perp}\cap\{y_1,y_2,\ldots
,y_n\}|+ \\ \nonumber & & (-1)\left
(n-|(w_1,w_2,\ldots,w_{m})^{\perp}\cap\{y_1,y_2,\ldots
,y_n\}|\right )\\ \nonumber
                        &=& q|(w_1,w_2,\ldots,w_{m})^{\perp}\cap \{y_1,y_2,\ldots ,y_n\}|-n,\\
\end{eqnarray*}
where $(w_1,w_2,\ldots,w_{m})^{\perp}=\{\langle
(x_1,x_2,\ldots,x_{m})\rangle\in \PG(m-1,q)\mid \sum_{i=1}^{m}
w_ix_i=0\}$. Using this formula for $\chi_{_{(w_1,w_2,\ldots,
w_{m})}}(\Omega)$ and Lemma~\ref{charsum-PDS}, one can prove the
following lemma.

\begin{lemma}
\label{geopartialset} Let ${\mathcal O}$ and $\Omega$ be defined
as above.  Then ${\mathcal O}$ is a projective $(n,m,h_1,h_2)$ set
in $\PG(m-1,q)$ if and only if $\chi_{_{(w_1,w_2,\ldots,
w_{m})}}(\Omega)=qh_1-n$ or $qh_2-n$, for every nontrivial
additive character $\chi_{_{(w_1,w_2,\ldots, w_{m})}}$,
$(w_1,w_2,\ldots, w_{m}) \in \Ff_q^{m}$. In other words,
${\mathcal O}$ is a projective $(n,m,h_1,h_2)$ set in $\PG(m-1,q)$
if and only if $\Omega$ is a $(q^{m}, (q-1)n, \lambda, \mu)$
partial difference set in the elementary abelian group
$(\Ff_q^{m},+)$, where
$\lambda=(q-1)n+(qh_1-n)(qh_2-n)+q(h_1+h_2)-2n$, and
$\mu=(q-1)n+(qh_1-n)(qh_2-n)$.
\end{lemma}

\subsection{Quadratic forms}
\label{quadraticforms}

Let $\Ff_q$ be the field of $q$ elements, where $q$ is a prime
power, and let $V$ be an $m$-dimensional vector space over
$\Ff_q$. A function $Q: V \rightarrow \Ff_q$ is called a
{\it{quadratic form}} if

\noindent(i). $Q(\alpha v) = \alpha^2 Q(v)$ for all $\alpha \in
\Ff_q$ and
$v \in V$, \\
(ii).  the function $B: V \times V \rightarrow \Ff_q$ defined by
$B(v_1,v_2) = Q(v_1+v_2) - Q(v_1) - Q(v_2)$ is bilinear.

We call $Q$ {\em{nonsingular}} if the subspace $W$ with the
property that $Q$ vanishes on $W$ and $B(w,v)=0$ for all $v\in V$
and $w\in W$ is the zero subspace. If the field $\Ff_q$ has odd
characteristic, then $Q$ is nonsingular if and only if $B$ is
nondegenerate; but this may not be true when $\Ff_q$ has
characteristic 2, because in that case $Q$ may not be zero on the
radical ${\rm Rad}(V)=\{w\in V\mid B(w,v)=0$ for all $v\in V\}$.
However, if $V$ is an even-dimensional vector space over an
even-characteristic field $\Ff_q$, then $Q$ is nonsingular if and
only if $B$ is nondegenerate (cf. \cite[p.~14]{cameron}).

Let $Q$ be a nonsingular quadratic form on an $m$-dimensional
vector space $V$ over $\Ff_q$. If $m$ is odd, then $Q$ is
equivalent to a quadratic form $x_1 x_2 + x_3 x_4 + \cdots +
x_{m-2} x_{m-1} + c x_{m} ^2$ for some scalar $c\in \Ff_q$ and it
is called a {\em{parabolic}} quadratic form.  If $m$ is even, then
$Q$ is equivalent to either a quadratic form $x_1x_2 + x_3 x_4 +
\cdots + x_{m-1} x_{m}$ (called {\em{hyperbolic}}, or type $+1$)
or $x_1 x_2 + x_3 x_4 + \cdots + x_{m-3}x_{m-2} +
p(x_{m-1},x_{m})$, where $p(x_{m-1},x_{m})$ is an irreducible
quadratic form in two indeterminates (called {\em{elliptic}}, or
type $-1$).

The {\it quadric} of the projective space $\PG(m-1,q)$
corresponding to a quadratic form $Q$ is the point set ${\mathcal
Q}=\{\langle v\rangle\in \PG(m-1,q)\mid Q(v)=0\}$. The following
theorems about the intersections of hyperplanes with quadrics in
$\PG(m-1,q)$ are well-known, see for example \cite{games},
\cite[p.~151]{br}, \cite{calderbankkantor}.

\begin{thm}  \label{ellquadform}
Let ${\mathcal Q}$ be a nonsingular elliptic quadric in
$\PG(2\ell-1,q)$. Then the hyperplanes of $\PG(2\ell-1,q)$
intersect ${\mathcal Q}$ in sets of two sizes $a$ and $b$, where
$$a=1+\frac {q(q^{\ell-1}+1)(q^{\ell-2}-1)} {q-1},\; b=\frac {q^{2\ell-2}-1} {q-1}.$$
If we use $\Omega$ to denote the set of nonzero vectors in
$\Ff_q^{2\ell}$ corresponding to ${\mathcal Q}$, then for any
nontrivial additive character $\chi$ of  $\Ff_q^{2\ell}$, we have
$\chi(\Omega)=(q^{\ell-1}-1)-q^{\ell}$ or $(q^{\ell-1}-1)$
according as the hyperplane of $\PG(2\ell-1,q)$ corresponding to
$\chi$ meets ${\mathcal Q}$ in $a$ or $b$ points. That is,
$\Omega$ is a
$(q^{2\ell},(q^{\ell}+1)(q^{\ell-1}-1),q^{2\ell-2}-q^{\ell-1}(q-1)-2,q^{2\ell-2}-q^{\ell-1})$-negative
Latin square type PDS in the additive group of $\Ff_q^{2\ell}$.
\end{thm}

\noindent{\bf Remarks} (1). Let ${\mathcal H}$ be a hyperplane of
a projective space $\PG(m-1,q)$, and denote by $W$ a point outside
${\mathcal H}$. If ${\mathcal Q}$ is a nonsingular quadric in
${\mathcal H}$ then the set
$${\mathcal C}=\cup_{X\in {\mathcal Q}} (WX)$$
is called a {\it cone} with vertex $W$ over ${\mathcal Q}$. Here
$WX$ means the set of points on the line through $W$ and $X$.

(2). The hyperplanes meeting a nonsingular elliptic quadric
${\mathcal Q}$ in $\PG(2\ell-1,q)$ in sets of size $a$ are called
{\it tangent} hyperplanes; such a hyperplane meets ${\mathcal Q}$
in a cone over a nonsingular elliptic quadric in $\PG(2\ell-3,q)$.
Any nontangent hyperplane meets ${\mathcal Q}$ in a nonsingular
parabolic quadric in that hyperplane.

(3). Similarly, let $\Omega$ be the set of nonzero vectors in
$\Ff_q^{2\ell}$ corresponding to a nonsingular hyperbolic quadric
in $\PG(2\ell-1,q)$. Then for any nontrivial additive character
$\chi$ of  $\Ff_q^{2\ell}$, we have
$\chi(\Omega)=q^{\ell}-(q^{\ell-1}+1)$ or $-(q^{\ell-1}+1)$. That
is, $\Omega$ is a $(q^{2\ell},(q^{\ell}-1)(q^{\ell-1}+1),
q^{2\ell-2}+q^{\ell-1}(q-1)-2,q^{2\ell-2}+q^{\ell-1})$-Latin
square type PDS in the additive group of $\Ff_q^{2\ell}$.
\vspace{0.1in}

Even though parabolic quadrics do not give rise to PDSs, we will
also need their intersection patterns with hyperplanes.

\begin{thm}  \label{paraquadform}
Let ${\mathcal Q'}$ be a nonsingular (parabolic) quadric in
$\PG(2\ell-2,q)$ (so the size of ${\mathcal Q'}$ is $\frac
{q^{2\ell-2}-1} {q-1}$). Then the hyperplanes of $\PG(2\ell-2,q)$
intersects ${\mathcal Q'}$ in sets of three sizes $t$, $h$, and
$e$ with respective multiplicities $T$, $H$ and $E$, where
\begin{eqnarray*}
t=\frac {q^{2\ell-3}-1} {q-1}, &  T=\frac {q^{2\ell-2}-1} {q-1},\\
e =t-q^{\ell-2}, & E=\frac {q^{2\ell-2}-q^{\ell-1}} {2},\\
h =t+q^{\ell-2}, & H=\frac {q^{2\ell-2}+q^{\ell-1}} {2}.\\
\end{eqnarray*}
\end{thm}

\noindent{\bf Remarks} (1). The $T$ hyerplanes with intersection
size $t$ are called {\it tangent} hyperplanes to ${\mathcal Q}'$,
such a hyperplane meets ${\mathcal Q}'$ in a cone over a parabolic
quadric. Each of the $E$ hyperplanes with intersection size $e$
meets ${\mathcal Q}'$ in a nonsingular elliptic quadric in that
hyperplane, and each of the $H$ hyperplanes with intersection size
$h$ meets ${\mathcal Q}'$ in a nonsingular hyperbolic quadric in
that hyperplane.

(2). If we use $\Omega'$ to denote the set of nonzero vectors in
$\Ff_q^{2\ell-1}$ corresponding to ${\mathcal Q}'$, then for any
nontrivial additive character $\chi$ of  $\Ff_q^{2\ell-1}$, we
have $\chi(\Omega')=-1$, $-1-q^{\ell-1}$, or $-1+q^{\ell-1}$.
(Here the character values $-1$, $-1-q^{\ell-1}$, and
$-1+q^{\ell-1}$ correspond respectively to cone section, elliptic
section and hyperbolic section.) \vspace{0.1in}

As in the study of all other types of difference sets, one of the
central problems in the study of partial difference sets is that
for a given parameter set, which groups of the appropriate order
contain a partial difference set with these parameters.  As far as
we know, no examples are known of negative Latin square type PDSs
in nonelementary abelian groups, and this paper will construct the
first such PDSs. We will focus on the case $q=4$, with $\Ff_4=\{
0, 1, \alpha, \alpha^2 = \alpha+1 \}$. We define a quadratic form
on $\Ff_4^{2\ell}$
\begin{eqnarray*}
Q_{ \ell,j}(x_1, x_2, \ldots, x_{2 \ell}) &=& (\alpha
x_1^2 + x_1 x_2 + x_2^2) + (\alpha x_3^2 + x_3 x_4 + x_4^2 )+ \cdots \\
&{}& + (\alpha x_{2j-1}^2 + x_{2j-1} x_{2j} + x_{2j}^2) + x_{2j+1}
x_{2j+2} + x_{2j+3} x_{2j+4} + \cdots + x_{2 \ell-1} x_{2 \ell}.\\
\end{eqnarray*}

\begin{lemma} \label{Qtjresults}
Let $Q_{ \ell,j}(x_1, x_2, \ldots, x_{2\ell})$ be the quadratic
form on $\Ff_4^{2\ell}$ defined above. Then

\noindent (a). When $j \geq 2$, $Q_{\ell,j}$ is projectively
equivalent to $Q_{\ell,j-2}$. \\
(b). When $j$ is odd, $Q_{ \ell,j}$ is elliptic. When $j$ is even, $Q_{\ell,j}$ is hyperbolic.\\
(c). $Q_{\ell,j}$ is nonsingular.\\
(d). $Q_{ \ell,j}(0,x_2,x_3,\ldots, x_{2 \ell})$ is a nonsingular
parabolic quadratic form in $2 \ell-1$ indeterminates.\\
(e).  $Q_{ \ell,j}(x_1, x_2 , \ldots, x_{2i-1}, x_{2i-1} + x_{2i},
\ldots, x_{2j-1}, x_{2j}, x_{2j+1}, \ldots, x_{2 \ell})= Q_{
\ell,j}(x_1, x_2, \ldots ,x_{2 \ell})$ for any $1 \leq i \leq j$.

\end{lemma}

\bproof  (a). The mapping $x_{2j-3} \mapsto (x'_{2j-3} + x'_{2j-1}
+ x'_{2j}); x_{2j-2} \mapsto (\alpha x'_{2j-3} + x'_{2j-2});
x_{2j-1} \mapsto (x'_{2j-1} + x'_{2j}); x_{2j} \mapsto (\alpha
x'_{2j-3} + x'_{2j-2} + x'_{2j})$; all other $x_i \mapsto x'_i$ is
an invertible linear transformation from $Q_{\ell,j}$ to
$Q_{\ell,j-2}$ as required.

(b).  If $j$ is even, then $Q_{\ell,j}$ is projectively equivalent
to $Q_{\ell, 0} = x_1 x_2 + x_3 x_4 + \cdots + x_{2 \ell - 1} x_{2
\ell}$, which is hyperbolic.  Similarly, if $j$ is odd, then
$Q_{\ell,j}$ is projectively equivalent to $Q_{\ell, 1} = \alpha
x_1^2 + x_1 x_2 + x_2^2 + x_3 x_4 + x_5 x_6 + \cdots + x_{2 \ell -
1} x_{2 \ell}$, which is elliptic.

(c). Let $B(x,x')$ be the bilinear form associated with
$Q_{\ell,j}$. Straightforward computations show that
$$B(x,x')=x_1x_2'+x_2x_1'+\cdots +x_{2\ell -1}x_{2\ell}'+x_{2\ell}x_{2\ell -1}',$$
which is nondegenerate. Hence $Q_{\ell,j}$ is nonsingular.

(d). According to Theorem 22.2.1 in \cite{hthas}, we define the
matrix $A=[a_{ik}]$, where $a_{ii}=2a_i$, $a_{ki}=a_{ik}$ for
$i<k$. Here $a_1=1$, $a_2=\alpha$, $a_3=1$, $a_4=\alpha$, $a_5=1$,
$\ldots$, $a_{2j-2}=\alpha$, $a_{2j-1}=1$, $a_{2j}=\cdots
=a_{2\ell -1}=0$, and $a_{12}=a_{13}=a_{14}=\cdots =0$,
$a_{23}=1$, $a_{45}=1$, etc. View $A$ as a matrix over $\Zz$, and
view $\alpha$ as an indeterminate for the time being. Compute
$\Delta=\frac{1}{2}{\rm det}(A)=(4\alpha -1)^{j-1}(-1)^{\ell-j}$.
Now view $\Delta$ modulo 2, we have $\Delta\neq 0$, by part (i) of
Theorem 22.2.1 in \cite{hthas}, this shows that
$Q_{\ell,j}(0,x_2,x_3,\ldots ,x_{2\ell})$ is nonsingular, and it
is necessarily parabolic (note: the associated bilinear form for
this quadratic form is degenerate, so we need the more
sophisticated argument to demonstrate nonsingularity).

(e). Straightforward computation, which we omit.  \eproof

We note that part (e) of the previous lemma will be used in many
character sum computations to get a sum of $0$ over the pair of
elements $(x_1, x_2 , \ldots, x_{2i-1}, x_{2i}, \ldots, x_{2
\ell})$ and $(x_1, x_2 , \ldots, x_{2i-1}, x_{2i-1} + x_{2i},
\ldots, x_{2 \ell})$.

\subsection{Galois ring preliminaries}

We need to recall the basics of Galois rings.  Interested readers
are referred to Hammons, {\em et. al} \cite{hkcss} for more
details.  We will only use Galois rings over $\Z_4$. A {\em Galois
ring over $\Z_{4}$ of degree $t$}, $t\geq 2$, denoted $\GR(4,t)$,
is the quotient ring $\Z_{4}[x]/ \langle \Phi(x) \rangle$, where
$\Phi(x)$ is a basic primitive polynomial in $\Z_{4}[x]$ of degree
$t$. Hensel's lemma implies that such polynomials exist. If $\xi$
is a root of $\Phi(x)$ in $\GR(4,t)$, then $\GR(4,t)=\Z_4[\xi]$
and the multiplicative order of $\xi$ is $2^t-1$.  In this paper,
we will only need $\GR(4,2)$, and that has the basic primitive
polynomial $\Phi(x) = x^2+x+1$.

The ring $R=\GR(4,t)$ is a finite local ring with unique maximal
ideal $2R$, and $R/2R$ is isomorphic to the finite field
$\Ff_{2^t}$. If we denote the natural epimorphism from $R$ to
$R/2R\cong \Ff_{2^t}$ by $\pi$, then $\pi(\xi)$ is a primitive
element of $\Ff_{2^t}$.

The set ${\mathcal T} = \{ 0,1,\xi,\xi^2, \ldots, \xi^{2^t-2} \}$
is a complete set of coset representatives of $2R$ in $R$. This
set is usually called a {\em Teichm\"uller system} for $R$.  The
restriction of $\pi$ to ${\mathcal T}$ is a bijection from
${\mathcal T}$ to $\Ff_{2^t}$, and we refer to this bijection as
$\pi_{_{\mathcal T}}$. An arbitrary element $\beta$ of $R$ has a
unique 2-adic representation
$$\beta = \beta_1 + 2 \beta_2,$$
where $\beta_1, \beta_2 \in {\mathcal T}$.  Combining
$\pi_{_{\mathcal T}}$ with this 2-adic representation and
specializing to the case of $\GR(4,2)$, we get a bijection $F_k$
from $\Ff_4^{2 \ell}$ to $(\GR(4,2))^k \times \Ff_4^{2 \ell - 2k}$
defined by
\begin{eqnarray*}
F_k: (x_1, x_2, \ldots, x_{2k-1}, x_{2k}, \ldots, x_{2 \ell})
\mapsto & (\pi_{_{\mathcal T}}^{-1}(x_1) + 2 \pi_{_{\mathcal
T}}^{-1}(x_2), \ldots, \pi_{_{\mathcal T}}^{-1}(x_{2k-1}) +
2 \pi_{_{\mathcal T}}^{-1}(x_{2k}),\\
{}& x_{2k+1}, x_{2k+2},\ldots, x_{2 \ell}).\\
\end{eqnarray*}
The inverse of this map is the map $F_k^{-1}$ from $(\GR(4,2))^k
\times \Ff_4^{2 \ell - 2k}$ to $\Ff_4^{2 \ell}$,
\begin{eqnarray*}
F_k^{-1}: (\xi_1+2\xi_2,\ldots
,\xi_{2k-1}+2\xi_{2k},\xi_{2k+1},\ldots ,\xi_{2\ell})\mapsto
&(\pi_{_{\mathcal T}}(\xi_1),\pi_{_{\mathcal T}}(\xi_2),\ldots
,\pi_{_{\mathcal T}}(\xi_{2k-1}),
\pi_{_{\mathcal T}}(\xi_{2k}), \\
{}& \xi_{2k+1},\xi_{2k+2},\ldots ,\xi_{2\ell})\\
\end{eqnarray*}

To simplify the notation we will usually omit the subindex in the
bijection $\pi_{_{\mathcal T}}^{-1}:\Ff_{2^t}\rightarrow {\mathcal
T}$. (So from now on, $\pi^{-1}$ means the inverse of the
bijection $\pi_{_{\mathcal T}}:{\mathcal T}\rightarrow
\Ff_{2^t}$.) We will show in the next section how we can use $F_k$
as a character sum preserving bijection to construct a PDS in a
nonelementary abelian group, namely the additive group of
$(\GR(4,2))^k \times \Ff_4^{2 \ell - 2k}$.

The {\em Frobenius map} $f$ from $R$ to itself is the ring
automorphism $f: \beta_1+2\beta_2 \mapsto \beta_1^2 + 2
\beta_2^2$. This map is used to define the {\em trace} ${\rm Tr}$
from $R$ to $\Z_{4}$, namely, ${\rm Tr}(\beta) = \beta + \beta^f +
\cdots + \beta^{f^{t-1}}$, for $\beta\in R$.  We note here that
the Galois ring trace $\Tr: R\rightarrow \Zz_4$ is related to the
finite field trace $\tr:\Ff_{2^m}\rightarrow \Ff_2$ via
\begin{equation}\label{commutrace}
\tr\circ \pi=\pi\circ \Tr. \end{equation}
As a consequence, we have $\sqrt{-1}^{{\rm Tr}(2 x)} =
\sqrt{-1}^{2{\rm Tr}(x)} = (-1)^{{\rm Tr}(x)} = (-1)^{\pi \circ
{\rm Tr}(x)} = (-1)^{{\rm tr}(\pi(x))}$, for all $x\in R$. The
trace of a Galois ring can be used to define all of the additive
characters of the ring, as demonstrated in the following
well-known lemma.

\begin{lemma}
\label{tracegiveschar} Let $\psi$ be an additive character of $R$.
Then there is a $\beta \in R$ so that $\psi(x) = \sqrt{-1}^{{\rm
Tr}(\beta x)}$ for all $x \in R$.
\end{lemma}

Since we can write $\beta = \beta_1 + 2 \beta_2$ for $\beta \in
\GR(4,2)$, where $\beta_k \in {\mathcal T}, k=1,2$, we will use
the notation $\psi_{_{\beta}} = \psi_{_{\beta_1 + 2 \beta_2}}$
indicating the ring element used to define the character $x\mapsto
\sqrt{-1}^{{\rm Tr}(\beta x)}$. If $\beta_1 =0$ but $\beta_2\neq
0$, then $\psi_{_{2 \beta_2}}$ is a character of order 2 and
$\psi_{_{2 \beta_2}}$ is principal on $2R$. If $\beta_1 \neq 0$,
then $\psi_{_{\beta_1 + 2 \beta_2}}$ is a character of order 4 and
$\psi_{_{\beta_1 + 2 \beta_2}}$ is nonprincipal on $2R$.
Characters of $(\GR(4,2))^k \times \Ff_4^{2 \ell - 2k}$ will be
written
\begin{equation}\label{defchar}
\Psi_k = \psi_{_{(\beta_{1} + 2 \beta_{2}, \ldots, \beta_{2k-1}+ 2
\beta_{2k})}} \otimes \chi_{_{(w_{2k+1},w_{2k+2}, \ldots, w_{2
\ell})}}
\end{equation}
for $\beta_i \in {\mathcal T}, 1 \leq i \leq 2k$, and $w_{i} \in
\Ff_4, 2k+1 \leq i \leq 2 \ell$.

\section{Construction of Partial Difference Sets in Nonelementary Abelian 2-Groups}

We are now ready to state the main result of this paper.  For $0
\leq k \leq j \leq \ell, \ell \geq 1$, define
\begin{equation}
D_{\ell, j,k} = \{ F_k(x_1, x_2, \ldots, x_{2 \ell})\mid (x_1,
x_2, \ldots, x_{2 \ell}) \in \Ff_4^{2 \ell} \backslash \{ 0 \},
Q_{\ell, j}(x_1, x_2, \ldots, x_{2 \ell})=0 \}.
\end{equation}
That is,
\begin{eqnarray*}
D_{\ell,j,k}=&\{(\xi_1+2\xi_2,\ldots
,\xi_{2k-1}+2\xi_{2k},\xi_{2k+1},\ldots ,\xi_{2\ell})\in
(\GR(4,2)^k\times\Ff_4^{2\ell-2k})\setminus \{0\} \mid \\
&{} Q_{\ell,j}(F_k^{-1}(\xi_1+2\xi_2,\ldots
,\xi_{2k-1}+2\xi_{2k},\xi_{2k+1},\ldots ,\xi_{2\ell}))=0\}
\end{eqnarray*}

We can think of $D_{\ell,j,k}$ as a ``lifting'' of the set of
nonzero vectors  $(x_1,x_2,\ldots,x_{2 \ell})\in \Ff_4^{2\ell}$
satisfying $Q_{\ell, j}(x_1, x_2, \ldots, x_{2 \ell})=0$.

\begin{thm} \label{maintheorem}
For $j$ odd, $1 \leq k \leq j \leq \ell$, the set $D_{\ell, j,k}$
is a $(4^{2 \ell}, (4^{\ell}+1)(4^{\ell-1}-1), 4^{2 \ell-2} - 3
\cdot 4^{ \ell-1} - 2, 4^{2 \ell-2} - 4^{ \ell-1})$-PDS in
$\Z_4^{2k} \times \Z_2^{4 \ell - 4 k}$.  For $j$ even, $1 \leq k
\leq j \leq \ell$, the set $D_{\ell, j,k}$ is a $(4^{2 \ell},
(4^{\ell}-1)(4^{\ell -1}+1), 4^{2 \ell-2} + 3 \cdot 4^{ \ell-1}
 - 2, 4^{2 \ell-2} + 4^{ \ell-1})$-PDS in $\Z_4^{2k} \times
\Z_2^{4 \ell - 4 k}$.
\end{thm}

This theorem immediately leads to the following corollary, which
lists the first known negative Latin square type PDSs in
nonelementary abelian groups.

\begin{cor}
\label{maincorollary} There are $(4^{2 \ell},
(4^{\ell}+1)(4^{\ell-1}-1), 4^{2 \ell-2} - 3\cdot 4^{\ell-1} - 2,
4^{2 \ell-2} - 4^{ \ell-1})$-negative Latin square type PDS in
$\Z_4^{2k} \times \Z_2^{4 \ell - 4 k}$ for every $k \leq \ell$
except possibly $\Z_4^{2 \ell}$ for $\ell$ even.
\end{cor}

By Lemma~\ref{charsum-PDS}, in order to prove
Theorem~\ref{maintheorem}, we demonstrate that all of the
nonprincipal characters of $(\GR(4,2))^k\times \Ff_4^{2 \ell - 2
k}$ have a sum over $D_{\ell, j,k}$ of $-4^{\ell -1} - 1 \pm 2
\cdot 4^{\ell -1}$ for $j$ odd and $4^{\ell -1} - 1 \pm 2 \cdot
4^{\ell -1}$ for $j$ even. By Theorem~\ref{ellquadform} and Remark
3 following that theorem, the set $F_k^{-1}(D_{\ell, j,k}) =
D_{\ell,j,0}$ is a PDS in the elementary abelian group of order
$4^{2 \ell}$, so it will have character sums equal to those we are
expecting for $D_{\ell, j,k}$.  The following sequence of lemmas
will indicate a connection between the character sums over
$D_{\ell, j,k}$ in the additive group of $(\GR(4,2))^k \times
\Ff_4^{2 \ell - 2 k}$ and the character sums over
$F_k^{-1}(D_{\ell, j,k})$ in the additive group of $\Ff_4^{2
\ell}$.  We first consider the characters of order two of the
additive group of $(\GR(4,2))^k \times \Ff_4^{2 \ell - 2 k}$.

\begin{lemma} \label{charactersoforder2}
Let $\Psi_k$ be a character of $(\GR(4,2))^k \times \Ff_4^{2 \ell
- 2 j}$ defined by (\ref{defchar}) with $\beta_{2i-1} = 0$ for $1
\leq i \leq k$. Then $\Psi_k(D_{\ell,j,k}) =
\chi_{_{(\pi(\beta_{2}),0,\pi(\beta_{4}),0,
\ldots,\pi(\beta_{2k}),0,w_{2k+1},w_{2k+2}, \ldots, w_{2
\ell})}}(F_k^{-1}(D_{\ell,j,k}))$.
\end{lemma}

\bproof Let $(\xi_{1} + 2 \xi_{2}, \xi_{3}+2 \xi_{4}, \ldots,
\xi_{2k-1}+ 2 \xi_{2k},\alpha_{2k+1},\alpha_{2k+2}, \ldots,
\alpha_{2 \ell}) \in D_{\ell,j,k}$. The character value of this
element is

\begin{eqnarray*}
& & \Psi_k(\xi_{1} + 2 \xi_{2}, \xi_{3}+2 \xi_{4}, \ldots,
\xi_{2k-1}+ 2 \xi_{2k},\alpha_{2k+1},\alpha_{2k+2}, \ldots,
\alpha_{2 \ell}) \\ & = & \sqrt{-1}^{\Tr(\sum_{i=1}^k 2 \beta_{2i}
\xi_{2i-1})}
(-1)^{\tr(\sum_{i'=2k+1}^{2 \ell} w_{i'} \alpha_{i'})} \\
& = & (-1)^{\tr(\sum_{i=1}^k \pi(\beta_{2i}) \pi(\xi_{2i-1}) +
\sum_{i'=2k+1}^{2 \ell} w_{i'}
\alpha_{i'})} \\
& = & \chi_{_{(\pi(\beta_{2}),0,
\ldots,\pi(\beta_{2k}),0,w_{2k+1},w_{2k+2}, \ldots, w_{2
\ell})}}(\pi(\xi_{1}), \pi(\xi_{2}),\ldots, \pi(\xi_{2k}),
\alpha_{2k+1}, \ldots, \alpha_{2 \ell})
\end{eqnarray*}

The second equality uses the fact that $\sqrt{-1}^{Tr(2 \beta_{2i}
\xi_{2i-1})} = (-1)^{tr(\pi(\beta_{2i} \xi_{2i-1}))}$ as mentioned
earlier in the discussion on trace. (Here $\Tr$ is the trace from
$\GR(4,2)$ to $\Z_4$, and $\tr$ is the trace from $\Ff_4$ to
$\Ff_2$.) This proves the lemma. \eproof

Thus, the character sums $\Psi_k(D_{\ell, j,k})$ associated to
characters $\Psi_k$ of order two will have the correct sum. In
order to prove Theorem~\ref{maintheorem}, we only need compute
$\Psi_k(D_{\ell, j,k})$, where $\Psi_k$ has order four.

Our strategy for proving Theorem~\ref{maintheorem} goes as
follows. We will first prove Theorem~\ref{maintheorem} in the case
$k=1$, then prove the whole theorem by strong induction on $k$. We
start by computing the character sum $\Psi_1(D_{\ell,j,1})$, where
$\Psi_1=\psi_{\beta_1+2\beta_2}\otimes \chi_{(w_3,w_4,\ldots
,w_{2\ell})}$ is a character of order 4, namely $\beta_1\neq 0$.
We will need the following definitions. Let
\begin{equation}\label{defomg0}
\Omega_0 = \{(2\xi_2,\xi_3,\ldots ,\xi_{2\ell})\in \GR(4,2) \times
\Ff_4^{2\ell-2}\setminus\{0\}\mid
Q_{\ell,j}(F_1^{-1}(2\xi_2,\xi_3,\ldots ,\xi_{2\ell}))=0\},
\end{equation}
and let
\begin{equation}\label{defO0}
O_0  =  F_1^{-1}(\Omega_0) =  \{(0,\pi(\xi_2),\xi_3,\ldots
,\xi_{2\ell})\in \Ff_4^{2\ell}\setminus\{0\}\mid
Q_{\ell,j}(0,\pi(\xi_2),\xi_3,\ldots ,\xi_{2\ell})=0\}.
\end{equation}

We observe that $\Psi_1(\Omega_0) =
\chi_{_{(\pi(\beta_{2}),\pi(\beta_{1}),w_{3},w_{4}, \ldots, w_{2
\ell})}}(O_0)$.  The next lemma shows their common sum when
$\beta_1 \neq 0$.

\begin{lemma}
\label{importantsubsetlemma} Suppose
$\chi_{_{(\pi(\beta_{2}),\pi(\beta_{1}),w_{3},w_{4}, \ldots, w_{2
\ell})}}$ is a character of $\Ff_4^{2 \ell}$ with $\beta_1,
\beta_2\in {\mathcal T}$ and  $\beta_{1}\neq 0$, and let $O_0$ be
as above. Then
$$\chi_{_{(\pi(\beta_{2}),\pi(\beta_{1}),w_{3},w_{4}, \ldots, w_{2 \ell})}}(O_0) = -1 \pm 4^{\ell-1}.$$
\end{lemma}

\bproof For convenience, let
$$O_0'=\{(x_2,x_3,\ldots ,x_{2\ell})\in \Ff_4^{2\ell -1}\setminus\{0\}\mid Q_{\ell,j}(0,x_2,\ldots ,x_{2\ell})=0\}.$$

By the definition of $O_0$, we have

$$\chi_{_{(\pi(\beta_{2}),\pi(\beta_{1}),w_{3},w_{4}, \ldots, w_{2 \ell})}}(O_0)=
\chi_{(\pi(\beta_1),w_3,\ldots ,w_{2\ell})}'(O_0').$$

Since $Q_{\ell,j}(0,x_2,x_3,\ldots ,x_{2\ell})$ is a nonsingular
parabolic quadratic form in $2\ell-1$ variables (cf.
Lemma~\ref{Qtjresults} (d)), the corresponding quadric in
$\PG(2\ell-2,4)$ has three intersection sizes with the
hyperplanes, leading to three distinct character values of $-1 \pm
4^{\ell-1}$ or $-1$ (see Theorem~\ref{paraquadform}, and the
remarks following that theorem). To prove the lemma, we need to
show that the condition $\beta_1 \neq 0$ excludes all characters
that have a sum of $-1$ over $O_{0}'$. To do that, we observe that
there are $4^{2 \ell - 2}-1$ nonprincipal characters satisfying
$\pi(\beta_{1}) = 0$, which is the same number of characters that
have a sum of $-1$ over $O_0'$ since the number of tangent
hyperplanes to a nonsingular parabolic quadratic form in
$\PG(2\ell-2,4)$ is $\frac {4^{2\ell-2}-1} {4-1}$ (cf.
Theorem~\ref{paraquadform}). We will show that all of these
characters $\chi_{(0,w_3,\ldots ,w_{2\ell})}'$ do indeed have a
sum of $-1$ over $O_0'$, implying that
$\chi_{(\pi(\beta_1),w_3,\ldots ,w_{2\ell})}'(O_0')$,
$\pi(\beta_1)\neq 0$, are equal to $-1 \pm 4^{\ell}$.

For any nontrivial character $\chi_{(0,w_3,\ldots ,w_{2\ell})}'$,
we have

\begin{eqnarray*}
\chi_{(0, w_3, \ldots , w_{2\ell})}'(O_0')&=&\sum_{x_2^2+x_3x_4+\cdots +x_{2\ell-1}x_{2\ell}=0}(-1)^{\tr(w_3x_3+\cdots +w_{2\ell}x_{2\ell})}\\
&=&\sum_{(0,0,\ldots ,0)\neq (x_3,x_4,\ldots ,x_{2\ell})\in \Ff_{4}^{2\ell-2}}(-1)^{\tr(w_3x_3+\cdots +w_{2\ell}x_{2\ell})}\sum_{x_2^2=x_3x_4+\cdots +x_{2\ell-1}x_{2\ell}}1\\
\end{eqnarray*}

Note that the inner sum in the last summation actually has only
one term. So
$$\chi_{(0, w_3, \ldots , w_{2\ell})}'(O_0')=\sum_{(0,0,\ldots ,0)\neq (x_3,x_4,\ldots ,x_{2\ell})\in \Ff_{4}^{2\ell-2}}(-1)^{\tr(w_3x_3+\cdots +w_{2\ell}x_{2\ell})}=-1.$$
As $(w_3,\ldots ,w_{2\ell})$ runs through all nonzero
$2\ell-2$-tuples, we get $-1$ as the character sum of $O_0'$
$(4^{2\ell-2}-1)$ times, these account for {\it all} the tangent
hyperplanes. Hence the lemma follows. \eproof

The following lemma considers the case $j$ odd and
$\chi_{(\pi(\beta_2),\pi(\beta_1),w_3,\ldots ,w_{2\ell})}(O_0) =
-1 - 4^{\ell-1}$, keeping in mind that the character sum over the
entire set $F_1^{-1}(D_{\ell,j,1})$ is $-1 - 4^{\ell-1} \pm 2
\cdot 4^{\ell-1}$ since $F_1^{-1}(D_{\ell,j,1})$ with $j$ odd is a
negative Latin square type PDS in the additive group of
$\Ff_4^{2\ell}$ (cf. Theorem~\ref{ellquadform}).

\begin{lemma}
\label{plusorminuslemma} Let $j$ be odd and
$\chi_{(\pi(\beta_2),\pi(\beta_1),w_3,\ldots ,w_{2\ell})}$ be a
character of $\Ff_4^{2\ell}$ with $\beta_1,\beta_2\in {\mathcal
T}, \beta_1 \neq 0$. If
$\chi_{(\pi(\beta_2),\pi(\beta_1),w_3,\ldots ,w_{2\ell})}(O_0) =
-1 - 4^{\ell-1}$, then
$$\chi_{(\pi(\beta_2),\pi(\beta_1),w_3,\ldots
,w_{2\ell})}(F_1^{-1}(D_{\ell,j,1}) \backslash O_0) = \pm 2 \cdot
4^{\ell-1}.$$
\end{lemma}

\bproof Obvious from the comments before the lemma. \eproof

The analogous lemma for the $j$ even case is given as follows.

\begin{lemma} Let $j$ be even and
$\chi_{(\pi(\beta_2),\pi(\beta_1),w_3,\ldots ,w_{2\ell})}$ be a
character of $\Ff_4^{2\ell}$ with $\beta_1,\beta_2\in {\mathcal
T}, \beta_1 \neq 0$. If
$\chi_{(\pi(\beta_2),\pi(\beta_1),w_3,\ldots ,w_{2\ell})}(O_0) =
-1 + 4^{\ell-1}$, then
$$\chi_{(\pi(\beta_2),\pi(\beta_1),w_3,\ldots
,w_{2\ell})}(F_1^{-1}(D_{\ell,j,1}) \backslash O_0) = \pm 2 \cdot
4^{\ell-1}.$$
\end{lemma}

\bproof Note that when $j$ is even, the character sum over the
entire set $F_1^{-1}(D_{\ell,j,1})$ is $-1 + 4^{\ell-1} \pm 2
\cdot 4^{\ell-1}$ since $F_1^{-1}(D_{\ell,j,1})$ with $j$ even is
a Latin square type PDS in the additive group of $\Ff_4^{2\ell}$.
The conclusion of the lemma is obvious from this observation.
\eproof

We now consider the other case from
Lemma~\ref{importantsubsetlemma}, namely that
$\chi_{(\pi(\beta_2),\pi(\beta_1),w_3,\ldots ,w_{2\ell})}(O_0) =
-1 + 4^{\ell-1}$ in the case $j$ odd; and
$\chi_{(\pi(\beta_2),\pi(\beta_1),w_3,\ldots ,w_{2\ell})}(O_0) =
-1 - 4^{\ell-1}$ in the case $j$ even.

\begin{lemma}
\label{zerosumlemma}Let $j$ be odd and let
$\chi_{(\pi(\beta_2),\pi(\beta_1),w_3,\ldots ,w_{2\ell})}$ be a
character of $\Ff_4^{2\ell}$ with $\beta_1, \beta_2\in {\mathcal
T}$, and $\beta_1 \neq 0$. If
$\chi_{(\pi(\beta_2),\pi(\beta_1),w_3,\ldots ,w_{2\ell})}(O_0) =
-1 + 4^{\ell-1}$, then
$$\chi_{(\pi(\beta_2),\pi(\beta_1),w_3,\ldots ,w_{2\ell})}(F_1^{-1}(D_{\ell,j,1}) \backslash O_0) = 0.$$
\end{lemma}

\bproof Since $j$ is odd, $F_1^{-1}(D_{\ell,j,1})$ corresponds to
a nonsingular elliptic quadric in $\PG(2\ell-1,4)$, hence its
character values $\chi_{(\pi(\beta_2),\pi(\beta_1),w_3,\ldots
,w_{2\ell})}(F_1^{-1}(D_{\ell,j,1}))$ are $-1-4^{\ell-1}\pm 2\cdot
4^{\ell-1}$. Assume to the contrary that
$\chi_{(\pi(\beta_2),\pi(\beta_1),w_3,\ldots
,w_{2\ell})}(F_1^{-1}(D_{\ell,j,1}) \backslash O_0)\neq 0.$ By the
assumption that $\chi_{(\pi(\beta_2),\pi(\beta_1),w_3,\ldots
,w_{2\ell})}(O_0) = -1 + 4^{\ell-1}$, we have
$$\chi_{(\pi(\beta_2),\pi(\beta_1),w_3,\ldots ,w_{2\ell})}(F_1^{-1}(D_{\ell,j,1}))=-1-4^{\ell-1}-2\cdot 4^{\ell-1},$$
that is, the hyperplane ${\mathcal H} : \pi(\beta_{2}) x_1 +
\pi(\beta_{1}) x_2 + w_3x_3+  \cdots + w_{2 \ell} x_{2 \ell}=0$
meets ${\mathcal Q}_{\ell,j}$ in a cone ${\mathcal C}$ with vertex
$W$ over a nondegenerate elliptic quadric in $PG(2 \ell-3,4)$.
(Here ${\mathcal Q}_{\ell,j}$ is the elliptic quadric defined by
$Q_{\ell,j}$ in $\PG(2\ell-1,4)$.) We write ${\mathcal H}\cap
{\mathcal Q}_{\ell,j}={\mathcal C}$.

Let ${\mathcal H}_{1}$ denote the hyperplane $x_{1}=0$. We know
that ${\mathcal H}_{1}\cap {\mathcal Q}_{\ell,j}={\mathcal
Q}'_{\ell,j}$ is a nonsingular parabolic quadric with equation
$x_2^2 +\alpha x_3^2+x_3x_4+x_4^2+\cdots +\alpha
x_{2j-1}^2+x_{2j-1}x_{2j}+x_{2j}^2+x_{2j+1}x_{2j+2}+ \cdots + x_{2
\ell - 1} x_{2 \ell}$ (cf. Lemma~\ref{Qtjresults} (d)).

Now consider ${\mathcal H}\cap {\mathcal Q}'_{\ell,j}$ (this will
determine the character value
$\chi_{(\pi(\beta_2),\pi(\beta_1),w_3,w_4,\ldots
,w_{2\ell})}(O_0)$). We have
$${\mathcal H}\cap {\mathcal Q}'_{\ell,j}=
{\mathcal H}_{1}\cap ({\mathcal H}\cap {\mathcal Q}_{\ell,j})=
{\mathcal H}_{1}\cap {\mathcal C}.$$

Note that ${\mathcal C}$ is a cone over an elliptic quadric, so
any hyperplane not through the vertex $W$ meets ${\mathcal C}$ in
an elliptic quadric \cite[p.~177]{br}. Hence if we can show that
${\mathcal H}_{1}$ does not go through the vertex $W$, then we
know that ${\mathcal H}_{1}\cap {\mathcal C}={\mathcal H}\cap
{\mathcal Q}'_{\ell,j}$ is not hyperbolic, thus
$\chi_{(\pi(\beta_2),\pi(\beta_1),w_3,w_4,\ldots
,w_{2\ell})}(O_0)\neq -1+4^{\ell-1}$, which is a contradiction.

The cone ${\mathcal C}$ is a degenerate quadric with a
1-dimensional radical being its vertex $W$. So we compute the
radical of ${\mathcal C}={\mathcal H}\cap {\mathcal Q}_{\ell,j}$.
Let $B(X,X')=Q_{\ell,j}(X+X')-Q_{\ell,j}(X)-Q_{\ell,j}(X')$ be the
bilinear form associated with $Q_{\ell,j}$, where $X=(x_1, x_2,
\ldots, x_{2 \ell})$ and $X'=(x'_1, x'_2, \ldots, x'_{2 \ell})$.
Then

\begin{eqnarray*}
{\rm Rad}({\mathcal H}\cap {\mathcal Q}_{\ell,j})&=&\{X \mid B(X,X')=0\; {\rm for}\; {\rm all}\; X'\in {\mathcal H}\}\\
&=&\{(x_1,x_2,\ldots, x_{2 \ell})\mid x_1 x'_2 + x_2 x'_1+\cdots +x_{2 \ell-1}x'_{2\ell}+x_{2 \ell}x'_{2 \ell - 1}=0, \\
& & {\rm for}\; {\rm all}\; (x'_1,x'_2,\ldots ,x'_{2 \ell})\in {\mathcal H}\}\\
&=&\{\epsilon ((\pi(\beta_{1}),\pi(\beta_2),w_{4},w_{3}, \ldots, w_{2 \ell}, w_{2 \ell - 1})\mid \epsilon\in \Ff_4^*\}\\
\end{eqnarray*}
Therefore the vertex of ${\mathcal C}$ is
$W=(\pi(\beta_{1}),\pi(\beta_{2}),w_{4},w_{3}, \ldots,
w_{2\ell},w_{2\ell-1})$. Since $\beta_{1}\neq 0$, we see that
${\mathcal H}_{1}: x_{1}=0$ does not go through $W$. The proof is
now complete. \eproof

The analogous lemma for the $j$ even case is given below.

\begin{lemma} Let $j$ be even and let
$\chi_{(\pi(\beta_2),\pi(\beta_1),w_3,\ldots ,w_{2\ell})}$ be a
character of $\Ff_4^{2\ell}$ with $\beta_1, \beta_2\in {\mathcal
T}$, and $\beta_1 \neq 0$. If
$\chi_{(\pi(\beta_2),\pi(\beta_1),w_3,\ldots ,w_{2\ell})}(O_0) =
-1 - 4^{\ell-1}$, then
$$\chi_{(\pi(\beta_2),\pi(\beta_1),w_3,\ldots ,w_{2\ell})}(F_1^{-1}(D_{\ell,j,1}) \backslash O_0) = 0.$$
\end{lemma}

The proof of this lemma is completely parallel to that of
Lemma~\ref{zerosumlemma}. We omit the details.

The above five lemmas are enough to prove
Theorem~\ref{maintheorem} in the case $k=1$ for both $j$ odd and
even (see the proof below). For the purpose of doing induction on
$k$, we define, for any integers $j,k$, $2\leq k < j \leq \ell$,
the following sets
\begin{eqnarray*}\label{defomgk}
\Upsilon_{k-1,k} & = &
\{(\xi_1+2\xi_2,\ldots,0+2\xi_{2k-2},0+2\xi_{2k},\xi_{2k+1},\ldots,\xi_{2\ell})\in
D_{\ell,j,k} \}; \\
\Upsilon_{k-1} & = &
\{(\xi_1+2\xi_2,\ldots,0+2\xi_{2k-2},\xi_{2k-1}+2\xi_{2k},\xi_{2k+1},\ldots,\xi_{2\ell})\in
D_{\ell,j,k}, \xi_{2k-1} \neq 0 \}; \\
\Upsilon_{k} & = &
\{(\xi_1+2\xi_2,\ldots,\xi_{2k-3}+2\xi_{2k-2},0+2\xi_{2k},\xi_{2k+1},\ldots,\xi_{2\ell})\in
D_{\ell,j,k}, \xi_{2k-3} \neq 0 \}; \\
\Upsilon  & = &
\{(\xi_1+2\xi_2,\ldots,\xi_{2k-3}+2\xi_{2k-2},\xi_{2k-1}+2\xi_{2k},\xi_{2k+1},\ldots,\xi_{2\ell})\in
D_{\ell,j,k}, \xi_{2k-1} \neq 0, \xi_{2k-3} \neq 0 \}; \\
U_{k-1,k}
 & = & \{(\xi_1+2\xi_2,\ldots,0,\pi(\xi_{2k-2}),0,\pi(\xi_{2k}),\xi_{2k+1},\ldots,\xi_{2\ell})\in
D_{\ell,j,k-2} \}; \\
U_{k-1} & = &
\{(\xi_1+2\xi_2,\ldots,0,\pi(\xi_{2k-2}),\pi(\xi_{2k-1}),\pi(\xi_{2k}),\xi_{2k+1},\ldots,\xi_{2\ell})\in
D_{\ell,j,k-2}, \pi(\xi_{2k-1}) \neq 0 \}; \\
U_{k} & = &
\{(\xi_1+2\xi_2,\ldots,\pi(\xi_{2k-3}),\pi(\xi_{2k-2}),0,\pi(\xi_{2k}),\xi_{2k+1},\ldots,\xi_{2\ell})\in
D_{\ell,j,k-2}, \pi(\xi_{2k-3}) \neq 0 \}; \\
U  & = &
\{(\xi_1+2\xi_2,\ldots,\pi(\xi_{2k-3}),\pi(\xi_{2k-2}),\pi(\xi_{2k-1}),\pi(\xi_{2k}),\xi_{2k+1},\ldots,\xi_{2\ell})\in
D_{\ell,j,k-2}, \\
& & \pi(\xi_{2k-1}) \neq 0, \pi(\xi_{2k-3}) \neq 0 \}; \\
\end{eqnarray*}

We observe that $D_{\ell, j, k} = \Upsilon_{k-1,k} \cup
\Upsilon_{k-1} \cup \Upsilon_{k} \cup \Upsilon$ and $D_{\ell, j,
k-2} = U_{k-1,k} \cup U_{k-1} \cup U_{k} \cup U$.  The following
lemma connects the character sum over $D_{\ell, j, k}$ with the
character sum over $D_{\ell, j, k-2}$.

\begin{lemma}
\label{k-2step} For $k \geq 2$, let $\Psi_k$ be defined as in
(\ref{defchar}), and let
$\Psi_{k-2}=\psi_{(\beta_1+2\beta_2,\ldots,\beta_{2k-5}+2\beta_{2k-4})}\otimes
\chi_{(\pi(\beta_{2k-2}),\pi(\beta_{2k-3}),\pi(\beta_{2k}),\pi(\beta_{2k-1}),w_{2k+1},\ldots,w_{2\ell})}$.
If $\beta_{2k-3}$ and $\beta_{2k-1}$ are both nonzero, then
$\Psi_k(\Upsilon_{k-1,k}) = \Psi_{k-2}(U_{k-1,k});
\Psi_k(\Upsilon_{k-1}) = - \Psi_{k-2}(U_{k-1});
\Psi_k(\Upsilon_{k}) = - \Psi_{k-2}(U_{k});$ and $\Psi_k(\Upsilon)
= \Psi_{k-2}(U)$.
\end{lemma}

\bproof  It is straightforward to see that
$\Psi_k(\Upsilon_{k-1,k}) = \Psi_{k-2}(U_{k-1,k})$.  We will show
that $\Psi_k(\Upsilon_{k-1}) = -\Psi_{k-2}(U_{k-1})$: the other
arguments are similar. We use part (e) of Lemma~\ref{Qtjresults}
to organize our sum to take advantage of pairs of the form
$\Psi_k(\xi_1+2\xi_2, \ldots, \xi_{2k-1}+2 \xi_{2k}, \xi_{2k+1},
\ldots, \xi_{2 \ell}) + \Psi_k(\xi_1+2\xi_2, \ldots, \xi_{2k-1}+2
(\xi_{2k-1}+\xi_{2k}), \xi_{2k+1}, \ldots, \xi_{2 \ell}) =
\Psi_k(\xi_1+2\xi_2, \ldots, \xi_{2k-1}+2 \xi_{2k}, \xi_{2k+1},
\ldots, \xi_{2 \ell})(1+(-1)^{tr(\pi(\beta_{2k-1})
\pi(\xi_{2k-1}))})$.  This last term will be $0$ unless
$\pi(\xi_{2k-1}) = \pi(\beta_{2k-1}^{-1})$.  We note that we will
have the same term $(1+(-1)^{tr(\pi(\beta_{2k-1})
\pi(\xi_{2k-1}))})$ in the $\Psi_{k-2}(U_{k-1})$ sum. Thus, this
character sum, while originally over all elements of
$\Upsilon_{k-1}$ or $U_{k-1}$, will only be over elements with
$\pi(\xi_{2k-1}) = \pi(\beta_{2k-1}^{-1})$. We can factor out the
term $i^{Tr(\beta_{2k-1} \xi_{2k-1})} = i^{Tr(\beta_{2k-1}
\beta_{2k-1}^{-1})} = i^{Tr(1)} = -1$ from the
$\Psi_k(\Upsilon_{k-1})$ sum, and what is left will be
$\Psi_{k-2}(U_{k-1})$.  Thus, $\Psi_k(\Upsilon_{k-1}) =
-\Psi_{k-2}(U_{k-1})$. \eproof

We will use strong induction, and hence we will assume that the
character sum $\Psi_{k-2}(D_{\ell,j,k-2})$ will have the correct
values.  The following lemma shows how we can get a correct sum
for $\Psi_k(\Upsilon_{k-1})$ and $\Psi_k(\Upsilon_{k})$ when
$\beta_{2k-3} = \beta_{2k-1} \neq 0$.

\begin{lemma}
\label{zerooutnegativeterms} If $\beta_{2k-3} = \beta_{2k-1} \neq
0$, then $\Psi_k(\Upsilon_{k-1}) = \Psi_k(\Upsilon_{k}) = 0$.
\end{lemma}

\bproof Suppose $\beta_{2k-3} = \beta_{2k-1} \neq 0$.  We will
show $\Psi_k(\Upsilon_{k-1}) = 0$: the other case is similar.
Using the fact that we only have to sum over elements of
$\Upsilon_{k-1}$ with $\xi_{2k-1}  = \beta_{2k-1}^{-1}$, and the
fact that $\pi(\xi_{2k-2}) = \alpha^2 \pi(\xi_1) + \pi(\xi_1)^2
\pi(\xi_2)^2 + \pi(\xi_2) + \cdots + \pi(\xi_{2k-4}) + \alpha^2
\pi(\xi_{2k-1}) + \pi(\xi_{2k-1})^2 \pi(\xi_{2k})^2 +
\pi(\xi_{2k}) +
\cdots+\alpha^2\xi_{2j-1}+\xi_{2j-1}^2\xi_{2j}^2+\xi_{2j}+\xi_{2 j
+ 1}^2 \xi_{2 j+2}^2+\cdots + \xi_{2 \ell - 1}^2 \xi_{2 \ell}^2$
from the quadratic form, we get

\begin{eqnarray*}
\Psi_k(\Upsilon_{k-1}) & = & \sum_{\xi_{2k-3}=0}
i^{Tr(\sum_{i=1}^{k} \beta_{2i-1} \xi_{2i-1})}
(-1)^{tr(\sum_{i=1}^{k} (\pi(\beta_{2i-1} \xi_{2i}) +
\pi(\beta_{2i} \xi_{2i-1}))+ w_{2k+1} \xi_{2k+1} +
\cdots + w_{2 \ell} \xi_{2 \ell})} \\
& = &  \sum_{\xi_{2k-3}=0} i^{Tr(\sum_{i=1}^{k-2} \beta_{2i-1}
\xi_{2i-1} + \beta_{2k-1} \xi_{2k-1})} (-1)^{tr(\sum_{i=1}^{k-2}
(\pi(\beta_{2i-1} \xi_{2i}) + \pi(\beta_{2i} \xi_{2i-1})))} \\
& & (-1)^{tr(\pi(\beta_{2k-3}) \{ \alpha^2 \pi(\xi_1) +
\pi(\xi_1)^2 \pi(\xi_2)^2 + \pi(\xi_2) + \cdots + \pi(\xi_{2k-4})
+ \alpha^2 \pi(\xi_{2k-1}) + (\pi(\beta_{2k-1})^{-1})^2
\pi(\xi_{2k})^2 + \pi(\xi_{2k}) + \cdots + \xi_{2 \ell - 1}^2
\xi_{2 \ell}^2
\})} \\
 & & (-1)^{tr(\pi(\beta_{2k-2})(0) + \pi(\beta_{2k-1} \xi_{2k}) +
\pi(\beta_{2k} \beta_{2k-1}^{-1}) +w_{2k+1}\xi_{2k+1}+
\cdots + w_{2 \ell} \xi_{2 \ell})} \\
\end{eqnarray*}

This last sum ranges over arbitrary values for $\xi_i$ except for
$i=2k-3$ and $i=2k-1$, where the values are fixed.  We can
rearrange the sum so there is an inner sum over all possible
values of $\xi_{2k}$.  If we do that, this inner sum will be
$\Sigma_{_{\xi_{2k}}} (-1)^{tr((\pi(\beta_{2k-1}) +
\pi(\beta_{2k-3}) + \pi(\beta_{2k-3}^2 \beta_{2k-1}^2))
\pi(\xi_{2k}))}$.  Since $\beta_{2k-3} = \beta_{2k-1}$, this sum
reduces to $\Sigma_{\xi_{2k}} (-1)^{tr(\pi(\beta_{2k-1})
\pi(\xi_{2k}))}$, and this sum is $0$, proving the lemma. \eproof

From Lemmas~\ref{k-2step} and \ref{zerooutnegativeterms} we
conclude that $\Psi_k(D_{\ell,j,k}) = \Psi_{k-2}(D_{\ell,j,k-2})$.
Induction will tell us the value of $\Psi_{k-2}(D_{\ell,j,k-2})$,
taking care of this case.

The final piece we need to prove the main theorem is to compute
the value of $\Psi_k(\Upsilon_{k-1})$ when $\beta_{2k-3} \neq
\beta_{2k-1}$, where both $\beta_{2k-3}$ and $\beta_{2k-1}$ are
nonzero. We have the same inner sum as in the proof of the last
lemma, but the inner sum in this case, namely $\Sigma_{\xi_{2k}}
(-1)^{tr((\pi(\beta_{2k-1}) + \pi(\beta_{2k-3}) +
\pi(\beta_{2k-3}^2 \beta_{2k-1}^2)) \pi(\xi_{2k}))}$ will be $4$
since $\pi(\beta_{2k-1}) + \pi(\beta_{2k-3}) + \pi(\beta_{2k-3}^2
\beta_{2k-1}^2) = 0$ for all possible choices.  We can extend this
idea to each pair $(\xi_{2i-1},\xi_{2i}), i \leq k-2$ to get
``inner sums'' of
\newline $\Sigma_{\xi_{2i-1},\xi_{2i}} i^{Tr(\beta_{2i-1}
\xi_{2i-1})} (-1)^{tr((\pi(\beta_{2i-1} \xi_{2i}) + \pi(\beta_{2i}
\xi_{2i-1}) + \pi(\beta_{2k-3}) ( \alpha^2 \pi(\xi_{2i-1}) +
\pi(\xi_{2i-1}^2 \xi_{2i}^2) + \pi(\xi_{2i}))))}$. We are assuming
$\beta_{2i-1} \neq 0$, since otherwise we could combine
$\Psi_k(D_{\ell,j,k}) = \Psi_{k-1}(D_{\ell,j,k-1})$ with the
inductive hypothesis to get the correct character sum. (Note here
that $\Psi_{k-1} =
\psi_{(\beta_1+2\beta_2,\ldots,\beta_{2k-3}+2\beta_{2k-2})}\otimes
\chi_{(\pi(\beta_{2k}),\pi(\beta_{2k-1}), w_{2k+1},\ldots, w_{2\ell})}$.) 
We can use Lemma~\ref{Qtjresults} (e) to restrict the possible values for
$\xi_{2i-1}$, and then consider the sums over $\xi_{2i-1}=0$ and
$\xi_{2i-1} = \beta_{2i-1}^{-1}$ separately. These sums are
$\Sigma_{\xi_{2i-1}=0, \xi_{2i}} (-1)^{tr((\pi(\beta_{2i-1}) +
\pi(\beta_{2k-3})) \pi(\xi_{2i}))}$ and $-(-1)^{tr(\pi(\beta_{2i}
\beta_{2i-1}^{-1}) + \alpha^2 \pi(\beta_{2k-3}
\beta_{2i-1}^{-1}))} \Sigma_{\xi_{2i-1}=\beta_{2i-1}^{-1},
\xi_{2i}} (-1)^{tr((\pi(\beta_{2i-1}) + \pi(\beta_{2k-3}) +
\pi(\beta_{2i-1}^2 \beta_{2k-3}^2)) \pi(\xi_{2i}))}$.  If
$\beta_{2i-1} = \beta_{2k-3}$, then the first sum is 4 and the
second sum is 0; if $\beta_{2i-1} \neq \beta_{2k-3}$, then the
first sum is 0 and the second sum is $\pm 4$.  Thus, in either
case, the total sum over the pair $(\xi_{2i-1}, \xi_{2i})$ is $\pm
4$.  When $k < j$, the ``inner sum'' for pairs
$(\xi_{2i-1},\xi_{2i}), k+1 \leq i \leq j$, will be \newline
$\Sigma_{\xi_{2i-1},\xi_{2i}}(-1)^{tr((w_{2i-1} \xi_{2i} + w_{2i}
\xi_{2i-1} + \pi(\beta_{2k-3}) ( \alpha^2 \xi_{2i-1} +
\xi_{2i-1}^2 \xi_{2i}^2 + \xi_{2i})))}$.  This is the character
sum of a quadratic form, which is $\pm 4$ (see \cite[p.~341,
Exercise 6.30]{ln}).

For $i> j$, we have ``inner sums'' of the form
$\Sigma_{\xi_{2i-1}, \xi_{2i}} (-1)^{tr(\pi(\beta_{2k-3})^2
\xi_{2i-1} \xi_{2i} + w_{2i-1} \xi_{2i-1} + w_{2i} \xi_{2i})}$.
Again this is the character sum of a quadratic form, which is
$\pm4$. Thus, we have $\ell - 1$ ``inner sums'', each of which are
$\pm4$. This proves the following lemma.

\begin{lemma}
\label{Upsilon_knonzero} Suppose $\beta_{2i-1} \neq 0$ for $i \leq
k$ and $\beta_{2k-3} \neq \beta_{2k-1}$. Then
$\Psi_k(\Upsilon_{k-1}) = \pm 4^{\ell - 1}$ and
$\Psi_k(\Upsilon_{k}) = \pm 4^{\ell - 1}$.
\end{lemma}

The sum $\Psi_k(\Upsilon_{k-1}) + \Psi_k(\Upsilon_k)$ is either
$0$ or $\pm 2 \cdot 4^{\ell - 1}$.  If the sum is $0$, then we can
use the same idea as in the comments after
Lemma~\ref{zerooutnegativeterms} to show that
$\Psi_k(D_{\ell,j,k}) = \Psi_{k-2}(D_{\ell,j,k-2})$, and induction
will tell us that this has the correct sum.  The following lemma
finishes the computations we need to prove the main theorem.

\begin{lemma}
\label{threecharacterlemma} Let $\Psi_k$ be defined as in
(\ref{defchar}) and $\Psi_{k-2}$ as in Lemma~\ref{k-2step} with
$\beta_{2k-3}\neq 0$ and $\beta_{2k-1}\neq 0$, and suppose
$D_{\ell,j,k-2}$ is a PDS (with parameters depending on $j$ odd or
even). If $j$ is odd and $\Psi_k(\Upsilon_{k-1}) +
\Psi_k(\Upsilon_k) = \pm 2 \cdot 4^{\ell-1}$, then
$\Psi_k(\Upsilon_{k-1,k}) + \Psi_k(\Upsilon) = -1 - 4^{\ell-1}$.
If $j$ is even and $\Psi_k(\Upsilon_{k-1}) + \Psi_k(\Upsilon_k) =
\pm 2 \cdot 4^{\ell-1}$, then $\Psi_k(\Upsilon_{k-1,k}) +
\Psi_k(\Upsilon) = -1 + 4^{\ell-1}$.
\end{lemma}

\bproof By Lemma~\ref{k-2step}, $\Psi_{k-2}(U_{k-1,k}) +
\Psi_{k-2}(U) = \Psi_{k-2}(D_{\ell,j,k-2}) - \Psi_{k-2}(U_{k-1}) -
\Psi_{k-2}(U_k) = \Psi_{k-2}(D_{\ell,j,k-2}) +
\Psi_{k}(\Upsilon_{k-1}) + \Psi_{k}(\Upsilon_k)$.  If $j$ is odd
and $\Psi_k(\Upsilon_{k-1}) + \Psi_k(\Upsilon_k) = 2 \cdot
4^{\ell-1}$, then $\Psi_{k-2}(U_{k-1,k}) + \Psi_{k-2}(U) = (-1 -
4^{\ell-1} \pm 2 \cdot 4^{\ell - 1}) + 2 \cdot 4^{\ell - 1}$ (the
case $\Psi_k(\Upsilon_{k-1}) + \Psi_k(\Upsilon_k) = -2 \cdot
4^{\ell-1}$ is similar, as is the $j$ even case). We will show
that $\Psi_{k-2}(U_{k-1,k}) + \Psi_{k-2}(U) = -1 + 3 \cdot 4^{\ell
- 1}$ leads to a contradiction. Define a character $\Psi'_{k-2}$
to have $\beta'_i = \beta_i$ except

\begin{description}
\item[(i)] $\beta'_{2k-2}$ chosen to satisfy $tr(\pi(\beta'_{2k-2}
\beta_{2k-3}^{-1})) = 1 + tr(\pi(\beta_{2k-2} \beta_{2k-3}^{-1}))$
\item[(ii)] $\beta'_{2k}$ chosen to satisfy $tr(\pi(\beta'_{2k}
\beta_{2k-1}^{-1})) = 1 + tr(\pi(\beta_{2k} \beta_{2k-1}^{-1}))$.
\end{description}

Neither $\pi(\beta_{2k-2})$ or $\pi(\beta_{2k})$ appears in
$\Psi_{k-2}(U_{k-1,k})$, so $\Psi'_{k-2}(U_{k-1,k}) =
\Psi_{k-2}(U_{k-1,k})$. Both $\pi(\beta_{2k-2})$ and
$\pi(\beta_{2k})$ appear in $\Psi_{k-2}(U)$ in the term
$(-1)^{tr(\pi(\beta_{2k-2} \beta_{2k-3}^{-1}) + \pi(\beta_{2k}
\beta_{2k-1}^{-1}))}$.  (Here again we use the fact that if
$\pi(\xi_{2i-1}) \neq \pi(\beta_{2i-1}^{-1})$ then
Lemma~\ref{Qtjresults} (e) will allow us to pair elements to sum
to $0$.) The conditions on $\beta'_{2k-2}$ and $\beta'_{2k}$ imply
that $(-1)^{tr(\pi(\beta'_{2k-2} (\beta'_{2k-3})^{-1}) +
\pi(\beta'_{2k} (\beta'_{2k-1})^{-1}))} =
(-1)^{tr(\pi(\beta_{2k-2} \beta_{2k-3}^{-1}) + \pi(\beta_{2k}
\beta_{2k-1}^{-1}))}$, so $\Psi'_{k-2}(U) = \Psi_{k-2}(U)$. The
element $\beta_{2k}$ appears in $\Psi_{k-2}(U_{k-1})$ in the term
$(-1)^{tr(\pi(\beta_{2k} \beta_{2k-1}^{-1}))}$, but $\beta_{2k-2}$
does not appear in this sum, implying that $\Psi'_{k-2}(U_{k-1}) =
-\Psi_{k-2}(U_{k-1})$. Similarly, $\Psi'_{k-2}(U_{k}) =
-\Psi_{k-2}(U_{k})$. Hence,

\begin{eqnarray*}
\Psi'_{k-2}(D_{\ell,j,k-2}) & = & \Psi'_{k-2}(U_{k-1,k})  +
\Psi'_{k-2}(U) +
\Psi'_{k-2}(U_{k-1}) + \Psi'_{k-2}(U_{k}) \\
& = & \Psi_{k-2}(U_{k-1,k}) + \Psi_{k-2}(U) - \Psi_{k-2}(U_{k-1})
-
\Psi_{k-2}(U_{k})\\
& = & \Psi_{k-2}(U_{k-1,k}) + \Psi_{k-2}(U) +
\Psi_{k}(\Upsilon_{k-1}) + \Psi_{k}(\Upsilon_{k})  \end{eqnarray*}

If $\Psi_{k-2}(U_{k-1,k}) + \Psi_{k-2}(U) = -1 + 3 \cdot 4^{\ell -
1}$, then $\Psi'_{k-2}(D_{\ell,j,k-2}) = -1 + 5 \cdot 4^{\ell -
1}$. We assumed that $D_{\ell,j,k-2}$ was a PDS, but $-1 + 5 \cdot
4^{\ell - 1}$ is not a correct character sum for this set, so that
contradiction proves the lemma. \eproof

\vspace{0.1in}

\noindent{\bf{Proof of Theorem~\ref{maintheorem}:}}  We proceed by
induction on $k$. For the $k=1$ case,
Lemma~\ref{charactersoforder2} shows that the characters $\Psi_1$
of order 2 have the correct character sums.  So we only need to
consider characters
$\Psi_1=\psi_{\beta_1+2\beta_2}\otimes\chi_{(w_3,w_4,\ldots,w_{2\ell})}$
of order 4 (i.e., $\beta_1 \neq 0$). We will only give the
detailed arguments in the case $j$ odd using Lemma 2.5 and 2.7.
The case $j$ even is similar (using Lemma 2.6 and 2.8).

We break $\Psi_1(D_{\ell,j,1})$ into three sums,

\begin{eqnarray*}
\Psi_1(D_{\ell,j,1})&=&\sum_{(\xi_1+2\xi_2,\xi_3,\ldots,\xi_{2\ell})\in D_{\ell,j,1}}(\sqrt{-1})^{\Tr((\beta_1+2\beta_2)(\xi_1+2\xi_2))}(-1)^{\tr(\sum_{i=3}^{2\ell}w_i\xi_i)}\\
&=&\sum_{(\xi_1+2\xi_2,\xi_3,\ldots,\xi_{2\ell})\in D_{\ell,j,1}, \Tr(\beta_1\xi_1)=0}(-1)^{\tr(\pi(\beta_2\xi_1+\beta_1\xi_2))}(-1)^{\tr(\sum_{i=3}^{2\ell}w_i\xi_i)}\\
&{}&-\sum_{(\xi_1+2\xi_2,\xi_3,\ldots,\xi_{2\ell})\in D_{\ell,j,1}, \Tr(\beta_1\xi_1)=2}(-1)^{\tr(\pi(\beta_2\xi_1+\beta_1\xi_2))}(-1)^{\tr(\sum_{i=3}^{2\ell}w_i\xi_i)}\\
&{}&+\sum_{(\xi_1+2\xi_2,\xi_3,\ldots,\xi_{2\ell})\in D_{\ell,j,1}, \Tr(\beta_1\xi_1)={\rm odd}}(\sqrt{-1})^{\Tr(\beta_1\xi_1)}(-1)^{\Tr(\pi(\beta_2\xi_1+\beta_1\xi_2))}(-1)^{\tr(\sum_{i=3}^{2\ell}w_i\xi_i)}\\
\end{eqnarray*}

The third sum above is over elements of $D_{\ell,j,1}$ with the
property that $\Tr(\beta_1 \xi_1)$ is odd. In that case, the
observation following Lemma 1.5 implies that each such element has
a matching element so that the pair will combine for a character
sum of $0$, so the sum over all these elements is $0$. We note
that this is also true for the sum of
$\chi_{(\pi(\beta_2),\pi(\beta_1),w_3,\ldots,w_{2 \ell})}$ over
the elements of $F_1^{-1}(D_{\ell,j,1})$ with $\tr(\pi(\beta_1)
\pi(\xi_1)) = 1$.

For the first sum above, note that $\Tr(\beta_1\xi_1)=0$ implies
that $\xi_1=0$, we see that the first sum is over the set of
elements of $D_{\ell,j,1}$ with $\xi_1 = 0$. Hence
\begin{eqnarray*}
\sum_{(\xi_1+2\xi_2,\xi_3,\ldots,\xi_{2\ell})\in D_{\ell,j,1}}
(\sqrt{-1})^{\Tr((\beta_1+2\beta_2)(2\xi_2))}(-1)^{\tr(\sum_{i=3}^{2\ell}w_i\xi_i)}&=&\Psi_1(\Omega_0)\\
&=&\chi_{(\pi(\beta_2),\pi(\beta_1),w_3,\ldots,w_{2\ell})}(O_0),\\
\end{eqnarray*}
where $\Omega_0$ and $O_0$ are defined in (\ref{defomg0}) and
(\ref{defO0}), respectively. By Lemma~\ref{importantsubsetlemma},
the first sum is equal to $-1 \pm 4^{\ell-1}$.

The second sum above is over elements of $D_{\ell,j,1}$ satisfying
$\xi_1 \neq 0$ and $\Tr(\beta_1 \xi_1)$ is even.  Since both
$\beta_1$ and $\xi_1$ are in the Teichm\"uller set ${\mathcal T}$,
their product is in ${\mathcal T}$ as well.  The only nonzero
element of ${\mathcal T}$ with an even trace is $1$, and $\Tr(1) =
2$.  This implies that the second sum is
$$S=\sum_{(\beta_1^{-1}+2\xi_2,\xi_3,\ldots,\xi_{2\ell})\in D_{\ell,j,1}}
\chi_{(\pi(\beta_2),\pi(\beta_1),w_3,\ldots,w_{2
\ell})}(\pi(\xi_1),\pi(\xi_2),\xi_3,\ldots,\xi_{2\ell}).$$
Combining this with the observation after our analysis of the
third sum, we see that
$$S=\chi_{(\pi(\beta_2),\pi(\beta_1),w_3,\ldots ,w_{2\ell})}(F_1^{-1}(D_{\ell,j,1}) \backslash O_0).$$

Hence we finally have

\begin{equation}\label{maineqn1}
\Psi_1(D_{\ell,j,1})=
\chi_{(\pi(\beta_2),\pi(\beta_1),w_3,\ldots,w_{2\ell})} (O_0) -
\chi_{(\pi(\beta_2),\pi(\beta_1),w_3,\ldots
,w_{2\ell})}(F_1^{-1}(D_{\ell,j,1}) \backslash O_0).
\end{equation}

Thus, using Lemmas~\ref{importantsubsetlemma},
{~\ref{plusorminuslemma} and ~\ref{zerosumlemma} we get
$\Psi_1(D_{\ell,j,1}) = -1 - 4^{\ell-1} \pm 2 \cdot 4^{\ell-1}$ as
required. This completes the proof of the theorem in the case
where $k=1$.

\vspace{0.1in} For the inductive step, suppose that
$D_{\ell,j,k'}$ is a PDS for all $k-1 \geq k' \geq 0, \ell \geq j
\geq k > 1$, with the appropriate parameters depending on $j$ odd
or even. To show that $D_{\ell,j,k}$ is also a PDS, we compute the
character sums $\Psi_{k}(D_{\ell,j,k})$, where
$$\Psi_{k}=\psi_{(\beta_1+2\beta_2,\ldots,\beta_{2k-1}+2\beta_{2k})}\otimes\chi_{(w_{2k+1},\ldots,w_{2\ell})}.$$
As in the $k=1$ case, Lemma~\ref{charactersoforder2} shows that
the characters $\Psi_{k}$ of order 2 have the correct character
sum. For characters $\Psi_{k}$ of order 4, we consider three
cases.

The first case is when there is a $\beta_{2i-1} = 0$ for some $1
\leq i \leq k$.  We can assume without loss of generality that
$i=k$, permuting the coordinates if necessary. In this case,
$\Psi_{k}(D_{\ell,j,k}) = \Psi_{k-1}(D_{\ell,j,k-1})$, and the
inductive hypothesis implies that this sum is correct (note here
that $\Psi_{k-1} =
\psi_{(\beta_1+2\beta_2,\ldots,\beta_{2k-3}+2\beta_{2k-2})}\otimes\chi_{(\pi(\beta_{2k}),
0, w_{2k+1},\ldots,w_{2\ell})}$ since $\pi(\beta_{2k-1}) = 0$).

The second case is when none of the $\beta_{2i-1}$ are $0$ for $i
\leq k$ but there is a pair $\beta_{2i-1} = \beta_{2i'-1}$ for
$i,i' \leq k, i \neq i'$.  Again, without loss of generality we
can assume that $i=k-1$ and $i'=k$.  By applying
Lemmas~\ref{k-2step} and \ref{zerooutnegativeterms}, we get

\begin{eqnarray*}
\Psi_{k}(D_{\ell,j,k}) & = & \Psi_{k}(\Upsilon_{k-1,k}) +
\Psi_{k}(\Upsilon_{k-1}) +\Psi_{k}(\Upsilon_{k})
+\Psi_{k}(\Upsilon) \\
& = & \Psi_{k}(\Upsilon_{k-1,k}) +\Psi_{k}(\Upsilon) \\
& = & \Psi_{k-2}(U_{k-1,k}) +\Psi_{k-2}(U) \\
& = & \Psi_{k-2}(D_{\ell,j,k-2})
\end{eqnarray*}

Induction tells us that $\Psi_k(D_{\ell,j,k})$ has the correct
value.

The final case is when none of the $\beta_{2i-1}$ are $0$ for $i
\leq k$, and all of the $\beta_{2i-1}$ are distinct.  According to
Lemma~\ref{Upsilon_knonzero}, there are $3$ possibilities for
$\Psi_{k}(\Upsilon_{k-1}) + \Psi_{k}(\Upsilon_{k})$: $0$ or $\pm 2
\cdot 4^{\ell-1}$.  If the sum is $0$, then the remarks after
Lemma~\ref{zerooutnegativeterms} indicate that
$\Psi_k(D_{\ell,j,k}) = \Psi_{k-2}(D_{\ell,j,k-2})$, so induction
implies that the character sum has the correct value.  When
$\Psi_{k}(\Upsilon_{k-1}) + \Psi_{k}(\Upsilon_{k}) = \pm 2 \cdot
4^{\ell-1}$, Lemma~\ref{threecharacterlemma} shows that
$\Psi_{k}(\Upsilon_{k-1,k}) + \Psi_{k}(\Upsilon)= -1 - 4^{\ell-1}$
when $j$ is odd, so $\Psi_k(D_{\ell,j,k}) = -1 - 4^{\ell-1} \pm 2
\cdot 4^{\ell-1}$ as required (the $j$ even case is similar).

Thus, in all cases the character sum $\Psi_k(D_{\ell,j,k})$ is as
required, proving the theorem. \eproof

\noindent{\bf Remarks} (1). We have checked that lifting a
nonsingular elliptic quadric in $\PG(3,q)$, where $q\neq 4$, in
the same way as in Theorem~\ref{maintheorem} (namely using a
straightforward bijection from multiple copies of the field to the
Teichm\"uller set representation of the Galois ring elements),
will not produce a PDS in nonelementary abelian 2-group.

(2). Since $D_{2,1,0}$ is a PDS in $G_0=(\Ff_4^4,+)$ (here
$D_{2,1,0}$ corresponds to an elliptic quadric in $\PG(3,4)$), the
Cayley graph $(G_0, D_{2,1,0})$ is a strongly regular graph with
vertex set $G_0$. Using GAP, it is found that the full
automorphism group of this graph has order $2^{12}\cdot 3^2\cdot
5\cdot 17$. Similarly, since $D_{2,1,1}$ is a PDS in
$G_1=\Zz_4^2\times\Ff_4^2$, the Cayley graph $(G_1, D_{2,1,1})$ is
a strongly regular graph with vertex set $G_1$. Again by using
GAP, it is found that the full automorphism group of this graph
has order $2^{12}\cdot 3^2\cdot 5$. Even though the two strongly
regular Cayley graphs $(G_0, D_{2,1,0})$ and $(G_1, D_{2,1,1})$
have the same parameters, they have different full automorphism
groups, hence they are nonisomorphic. Based on further
computations of automorphism groups, we conjecture that the
strongly regular Cayley graphs arising from $D_{\ell,j,k}$, $k>0$,
are never isomorphic to the classical strongly regular Cayley
graphs arising from quadratic forms over $\Ff_4$.

(3) The PDSs $D_{\ell,j,k}$ and $D_{\ell,j',k}$ are equivalent if
$j-j'$ is even.  The mapping defined in the proof of
Lemma~\ref{Qtjresults} (a) induces a group automorphism $\phi$ on
$(\GR(4,2))^k\times \Ff_4^{2 \ell - 2 k}$ satisfying
$\phi(D_{\ell,j,k}) = D_{\ell,j-2,k}$ for $0 \leq k \leq j-2$
($\phi$ will fix all elements of the group $(\GR(4,2))^k\times
\Ff_4^{2 \ell - 2 k}$ except $x_{2j-3}$ through $x_{2j}$, and
these are mapped as directed in the proof).  We chose to do the
proof in the general form since that made the induction easier to
read, but we could have stated the result in terms of
$D_{\ell,j,j}$ and $D_{\ell,j+1,j}$ and gotten a PDS in each
equivalence class.

\section{Future Directions}

The following list describes possible implications of the results
in this paper.

(1).  Are there other groups of the appropriate order that support
PDSs with the same parameters in this paper?  In particular, are
there other combinations of $\Z_4$ and $\Z_2$ that contain PDSs? A
good candidate might be to see if $\Z_4^4$ contains a PDS with
negative Latin square parameters, finishing off the one case we
can't do in Corollary~\ref{maincorollary}.

(2).  Can the ``bijection idea'' produce any more sets where the
character sums are so well preserved?

(3).  Can we do anything in other contexts to get groups such as
$\Z_8$ involved in a PDS group?

 \vspace{0.5cm}

\noindent{\bf Acknowledgement:} We thank Frank Fiedler for helping
us compute the automorphism groups of some strongly regular
graphs. The research of Qing Xiang was supported in part by NSA
grant MDA 904-99-1-0012.  We would also like to thank the
anonymous referee for suggestions that improved the exposition of
this paper.

\end{document}